\documentclass{emsprocart}
\usepackage{mathrsfs,amsmath}
\usepackage[all]{xy}
\usepackage[pdftex,bookmarks]{hyperref}
\newcommand{\Z}{\mathbb{Z}}

\newcommand{\R}{\mathbb{R}}
\newcommand{\A}{\mathscr{A}}
\newcommand{\B}{\mathscr{B}}

\newcommand{\PP}{\mathcal{P}}
\newcommand{\E}{\mathcal{E}}
\newcommand{\ootimes}{\widehat{\otimes}}

\title[CLIFFORD MODULES AND TWISTED $K$-THEORY]{CLIFFORD MODULES AND TWISTED $K$-THEORY}
\author[Max KAROUBI]{by Max KAROUBI}
\contact[max.karoubi@gmail.com]
{Max Karoubi\\
Universit\'e Paris 7 - Denis Diderot\\
2 place Jussieu\\
75251 PARIS (FRANCE)}
\begin{document}
The purpose of this short paper is to make the link between the fundamental work of Atiyah, Bott and Shapiro \cite{ABS} and twisted $K$-theory as defined by P. Donovan, J. Rosenberg and the author \cite{DK} \cite{R} \cite{K5}. This link was implicit in the literature (for bundles over spheres as an example) but was not been explicitly defined before.\\ \\
The setting is the following: $V$ is a real vector bundle on a compact space  $X$, provided with a non degenerate quadratic form to which we associate a bundle of (real or complex) Clifford algebras denoted by $C(V)$; the quadratic form is implicit in this notation. We denote by $M(V)$ the Grothendieck group associated to the category of (real or complex) vector bundles provided with a structure of (twisted) $\Z/2$-graded $C(V)$-module. Another way to describe $M(V)$ is to consider the bundle $V\oplus
1$, where the symbol ``1" denotes the trivial vector bundle of rank one with a positive quadratic form. Then $M(V)$ is just the Grothendieck group  $K(\Lambda_1)$ of the category $\PP(\Lambda_1)$ whicho objects are finitely generated projective modules over $\Lambda_1$. The notation $\Lambda_n$ means in general $\Lambda \ootimes C^{0,n}$, where $\Lambda$ is the ring of continuous sections of the $\Z/2$-graded bundle $C(V)$ and $C^{0,n}$ is the Clifford algebra of $\R^n$ with a positive quadratic form.\\ \\
Following \cite{ABS}, we define $A(V)$ as the cokernel of the homomorphism induced by restriction of the scalars :
$$A(V) = Coker [M(V\oplus 1) \to M(V)]  = Coker[K(\Lambda_2) \to
K(\Lambda_1) ].$$

\subsection*{Remark} Let us denote by $V^-$ the vector bundle $V$ with the opposite quadratic form. It is quite easy to see\footnote{If $(v, t)$ is a symbol for the action of $V\oplus 1$, with $t^2 = 1$, we change it in $(vt, t)$ which represents the action of $V^-\oplus 1$,} that the category of $C(V\oplus 1)$-modules is isomorphic to the category of $C(V^-\oplus 1)$-modules. From now on, we assume that the quadratic form on $V$ is positive (in which case $A(V^-)$ was the original definition of \cite{ABS}).\\ \\
With these definitions, we have the following theorem, where $K(V)$ denotes the real or complex reduced $K$-theory of the Thom space of $V$.

\begin{theorem}{} There is an exact sequence between $K$-groups\footnote{One might also write $K_{-1}$ instead of $K^1$.}
$$
K(\Lambda_2) \to
K(\Lambda_1)  \to
K(V)  \to
K^1(\Lambda_2) \to 
K^1(\Lambda_1)$$

In particular, $A(V)$ is a subgroup of  $K(V)$ which coincides with it in the following important cases :
\begin{itemize}
\item[a)] $K^1(\Lambda_2) = 0$, for instance when $X$  is reduced to a point.
\item[b)] $V$ is oriented of rank divisible by 4.
\item[c)] $V$ is oriented of even rank in the framework of complex $K$-theory.
\end{itemize}
\end{theorem}
\begin{proof} According to the general theory developed in \cite{K1}, $K(V) \equiv
K^1(V\oplus 1)$ is canonically isomorphic to the $K^1$-group of the Banach functor

$$\phi : \E^{V\oplus 2} (X)  \to \E^{V\oplus 1} (X),$$
where $\E^W (X)$ denotes the category of vector bundles provided with a $C(W)$-module structure. According to the Serre-Swan theorem, the categories involved are equivalent to categories $\PP(R)$ for suitable rings $R$, in this case $\Lambda_2$ or $\Lambda_1$. The first part of the theorem follows from these general considerations.\\
If $X$  is a point, the category $\PP(R)$ is finite dimensional and therefore its $K^1$-group is trivial.
On the other hand, if $V$ is oriented of rank $n$ divisible by 4, let us choose an orthonormal \textsl{oriented} basis $e_1, \dots,  e_n$ on each fiber $V_x$, $x\in X$. Then the product $\varepsilon = e_1\dots e_n$ in the Clifford algebra $C(V_x)$ is independant of the choice of the basis since $\varepsilon$ commutes with the action of $SO(n)$ and defines therefore a continuous section of $C(V)$. On the other hand, for any $W$, there is an isomorphism between the graded tensor product $C(V)\ootimes C(W)$ and the nongraded one $C(V)\otimes C(W)$. In order to see it, we send $V\oplus W$ to $C(V)\otimes C(W)$ by the formula
$$(v, w) \mapsto v \otimes 1 + \varepsilon\otimes w,$$
The fact that $n$ is even shows that $\varepsilon$ anticommutes with $v$. Moreover, if 4 divides $n$, the square of $\varepsilon$ is 1. Therefore, by the universal property of Clifford algebras, the previous map induces the required isomorphism $C(V)\ootimes C(W) \equiv C(V) \otimes C(W)$. If $n = 4k + 2$ and in the framework of complex $K$-theory, one may replace $\varepsilon$ by  $\varepsilon\sqrt{-1}$ in order to get the same result.\\
In our situation, $W$ is of dimension one or two and the Banach functor
$$\PP(\Lambda)\sim
\PP(\Lambda\ootimes C^{0,2}) \to
\PP(\Lambda\ootimes C^{0,1})  \sim
\PP(\Lambda) \times \PP(\Lambda)
$$
may be identified to the diagonal functor through the previous category isomorphisms. This shows that the map
$$K^1(\Lambda \ootimes C^{0,2}) \to K^1(\Lambda \ootimes C^{0,1})$$
is injective and concludes the proof of the theorem.
\end{proof}
%%%
\subsection*{Example} When  $X$  is reduced to a point, the theorem implies that the reduced $K$-theory of the sphere $S^n$ is the cokernel of the map
$$K(C^{0,n+2}) \to K(C^{0,n+1})$$
which is the same as the cokernel of the map
$$K(C^{n+1,1}) \to K(C^{n,1})
$$
as we noticed earlier. This is the starting remark in \cite{ABS} which was the inspiration of \cite{K1}, where the notation $C^{p,q}$ is used.
%%%
\subsection*{Generalizations} Since the main tool used here is the real Thom isomorphism proved in \cite{K1} and \cite{K3}, the previous theorem might be generalized to the equivariant case. For instance, if $G$ is a finite group acting linearly on $\R^n$, the group $K_G(\R^n)$ is isomorphic to the cokernel of the following map
          $$K(G\ltimes C^{0,n+2}) \to K(G \ltimes C^{0,n+1})$$
where the involved rings are crossed products of $G$ by Clifford algebras. More precise results may be found in \cite{K4}.\\
\\
\\
Another generalization is to consider modules over bundles of $\Z/2$-graded Azumaya algebras $\A$ as in \cite{DK}, instead of bundles of Clifford algebras. The analog of the group $A(V)$ is now what we might call the ``algebraic twisted $K$-theory" of $\A $, denoted by  $K^{\A}_{alg}(X)$ which is the cokernel of the map $K(\Lambda_2) \to K(\Lambda_1)$, where $\Lambda$ denotes the ring of continuous sections of the $\Z/2$-graded algebra bundle $\A$. We can prove, as in the previous theorem, that this new group is a subgroup\footnote{One has to use again the Thom isomorphism in twisted $K$-theory as stated in \cite{K5}. } of the usual twisted $K$-theory of  $X$  denoted by $K^\A(X)$. It coincides with it in some important cases, for instance if $\A$  is oriented (in the graded sense) with fibers modelled on matrix algebras over the real or complex numbers.

\subsection*{The multiplicative structure}
It is well known that the twisted $K$-groups $K^\A(X)$ can be provided with a cup-product structure (see \cite{DK} \S  7 and also \cite{K2}). This cup-product makes a heavy use of Fredholm operators in Hilbert spaces. This machinery is unavoidable, especially for the odd $K$-groups.\\ \\
More precisely, as shown in \cite{DK}, the elements which define the group $K^\A(X)$  are pairs $(E, D)$ where $E$ is a $\Z/2$-graded Hilbert bundle provided by an $\A$-module structure and 
$D : E \to E$ is a family of Fredholm operators which are self-adjoint, of degree one and commute  (in the graded sense) with the action of $\A$. According to the Thom isomorphism in twisted $K$-theory \cite{K5}, $K^\A(X)$ is also the $K^1$-group of the Banach functor
$$\phi : \E^{\A\ootimes C^{0,2}}(X) \to\E^{\A\ootimes C^{0,1}}(X).$$
We have therefore the following exact sequence (as for Clifford modules)
$$K(\A\ootimes C^{0,2}) \to K(\A\ootimes C^{0,1})\to K^\A(X)\to K^1(\A\ootimes C^{0,2})\to
K^1(\A\ootimes C^{0,1}).$$

A closer look at the connecting homomorphism
$$K(\A\ootimes C^{0,1}) \to K^\A(X)$$
shows that it associates the pair $(E,0)$ to a (finite dimensional) vector bundle $E$ which is a module over $A\ootimes C^{0,1}$. In other words, the elements in  $K^\A(X)$ corresponding to finite dimensional bundles $E$ are just elements of the cokernel of the map $K(\A\ootimes C^{0,2}) \to K(\A\ootimes C^{0,1})$
which we might call $A(\A)$, if we follow the conventions of \cite{ABS} or  simply  $K^\A_{alg}(X)$ as we did before.\\ \\
Since the elements of  $K^\A_{alg}(X)$ are associated to finite dimensional bundles (with the Fredholm operator reduced to 0), the usual cup-product
$$K^\A(X) \times K^{\A'}(X) \to K^{\A\ootimes \A' }(X)$$
induces a pairing between the algebraic parts
$$K^\A(X)_{alg} \times K^{\A'}_{alg}(X) \to K^{\A\ootimes \A' }_{alg}(X).$$
On the other hand, it might be interesting to characterize the elements of  $K^\A(X)$  which are ``algebraic". They belong to the kernel of the map
$$\phi :  K^\A(X)\to K^1(\A\ootimes C^{0,2}) = K(\B/\A\ootimes C^{0,2}).$$

The notation $\B/\A$ represents here the bundle of Calkin algebras associated to $\A$ (the structural group of $\A$ is $PU(H)$ where $H$ is an infinite dimensional Hilbert space, as mentioned in \cite{K5} and \cite{R}). This map $\phi$ is easy to describe: it associates to a couple $(E, D)$ as before the space of sections of the $\B/\A-$bundle associated to $E$. It is provided with the action of $C^{0,2}$ described by the grading and the involution on the Calkin bundle induced by the polar decomposition of $D$. This description of the ``algebraic" elements also holds in the generalization of twisted $K$-theory considered by J. Rosenberg \cite{R}\cite{K5}.
%REFERENCES


\begin{thebibliography}{99}
\bibitem{ABS} M.F. ATIYAH, R. BOTT and A. SHAPIRO. Clifford modules. Topology 3, pp. 3-38 (1964).

\bibitem{DK} P. DONOVAN and M. KAROUBI. Graded Brauer groups and $K$-theory with local coefficients. Publ. Math. IHES 38, pp. 5-25 (1970). French summary in : Groupe de Brauer et coefficients locaux en $K$-th\'eorie. Comptes Rendus Acad. Sci. Paris, t. 269, pp. 387-389 (1969).

\bibitem{K1} M. KAROUBI. Alg\`ebres de Clifford et $K$-th\'eorie. Ann. Sci. Ecole Norm. Sup. (4), pp. 161-270 (1968).

\bibitem{K2} M. KAROUBI. Alg\`ebres de Clifford et op\'erateurs de Fredholm. Springer Lecture Notes in Maths N° 136, pp. 66-106 (1970).  Summary in Comptes Rendus Acad. Sci. Paris, t. 267, pp. 305 (1968).

\bibitem{K3} M. KAROUBI. Sur la $K$-th\'eorie \'equivariante. Springer Lecture Notes in Math. N° 136, pp. 187-253 (1970).

\bibitem{K4} M. KAROUBI. Equivariant $K$-theory of real vector spaces and real projective spaces. Topology and its applications 122, pp. 531-546 (2002).

\bibitem{K5} M. KAROUBI. Twisted $K$-theory, old and new. \href{http://arxiv.org/abs/math/0701789}{ArXiv math 0701789}  (to appear in the Journal of European Math. Society in 2008).

\bibitem{R} J. ROSENBERG. Continuous-trace algebras from the bundle theoretic point of view. J. Austral. Math. Soc. A 47, pp. 368-381 (1989).

 \end{thebibliography}
\end{document}